\documentclass[runningheads]{llncs}

\usepackage{mystyle}
\usepackage[T1]{fontenc}
\usepackage{graphicx,xcolor, subcaption, multirow, booktabs, threeparttable}

\begin{document}
\title{On the Convexification of a Class of Mixed-Integer Conic Sets}
%
%
\author{Guxin Du \and
Rui Chen \and
Linchuan Wei}
\authorrunning{G. Du et al.}
%
\institute{School of Data Science, The Chinese University of Hong Kong, Shenzhen
\email{guxindu@link.cuhk.edu.cn, \{rchen,linchuanwei\}@cuhk.edu.cn}}
\maketitle              
\begin{abstract}
We investigate mixed-integer second-order conic (SOC) sets with a nonlinear right-hand side in the SOC constraint, a structure frequently arising in mixed-integer quadratically constrained programming (MIQCP). Under mild assumptions, we show that the convex hull can be exactly described by replacing the right-hand side with its concave envelope. This characterization enables strong relaxations for MIQCPs via reformulations and cutting planes. Computational experiments on distributionally robust chance-constrained knapsack variants demonstrate the efficacy of our reformulation techniques.

\end{abstract}
%
%
%


\section{Introduction}\label{sec:intro}

In this paper, we consider a mixed-integer nonlinear set specified by a second-order conic (SOC) constraint with a nonlinear right-hand side. Specifically, we consider
\begin{equation}\label{eq:set_Z}
    \calZ := \left\{(\vecx,\vecy) \in \calX \times \mathbb{R}^m: \norm{A\vecx + B\vecy + \vecd}_2 \leq f(\vecx) \right\},
\end{equation}
where $\calX \subsetneq \bbR^n$ is a compact set, $A$ and $B$ are matrices of appropriate dimensions, $\vecd \in \bbR^p$ is a constant vector, and $f: \calX \to \bbR$ is a function that is upper semi-continuous on $\calX$. We are particularly interested in the case when $\vecx$ are binary variables, i.e., $\calX \subseteq\{0,1\}^n$.

This general $\calZ$ encompasses, as a special case, mixed-integer convex quadratically constrained sets, which form the core structure in mixed-integer chance-constrained and distributionally robust chance-constrained programming \cite{zhangAmbiguousChanceConstrainedBinary2018}. Such models have diverse applications, including surgery planning \cite{deng2019chance} and energy-efficient cloud computing \cite{shen2014stochastic}. For particular choices of $A$, $B$, $d$, and $f$, $\calZ$ reduces to the following mixed-binary second-order conic set studied in \cite{atamturk2020submodularity}:
\begin{equation*}
\left\{(\vecx,\vecy) \in  \{0,1\}^{n} \times \mathbb{R}^{m}: \sqrt{\sigma + \sum_{i = 1}^{n} a_i x_i + \sum_{i=1}^{m} b_i y_i^2} \leq T\right\},
\end{equation*}
whose convex hull description remains an open question. 

Our contribution lies within the broader literature on strong formulations for mixed-integer nonlinear sets, which can be classified into three streams. 
One stream of research extends techniques from the linear to the nonlinear setting, including Gomory cuts~\cite{ccezik2005cuts}, mixed-integer rounding cuts~\cite{atamturk2010conic}, lifting~\cite{atamturk2011lifting}, intersection cuts~\cite{modaresi2016intersection}, disjunctive programming~\cite{ceria1999convex}, and lift-and-project cuts~\cite{stubbs1999branch}. 
Another line of work exploits the submodularity~\cite{edmonds2003submodular} inherent in certain structures to derive tight relaxations~\cite{atamturk2020submodularity,atamturk2008polymatroids,atamturk2009submodular}. 
Most relevant to our study is the third stream, which focuses on strong formulations for \emph{structured} mixed-integer nonlinear sets. 
These include quadratically constrained sets~\cite{burer2017convexify,burer2020exact,modaresi2017convex,wang2022tightness}, quadratic programs with semicontinuous variables~\cite{akturk2009strong,anstreicher2021quadratic,de2024explicit,frangioni2006perspective,gunluk2010perspective,wei2024convex}, mixed-integer conic sets~\cite{atamturk2020submodularity,atamturk2010conic,gomez2021strong,kilincc2025conic,shafiee2024constrained}, and nonlinear disjunctive sets~\cite{kilincc2015two}, among others. 

In Section~\ref{sec:theory}, we prove that under some mild assumptions, the convex hull of the set $\calZ$ is given by  \
\begin{equation} \label{eq:set_W}
        \mathcal{W}: = \left\{(\vecx,\vecy) \in \conv(\calX) \times \bbR^m: \norm{A\vecx + B\vecy + \vecd}_2 \leq \hat{f}(\vecx) \right\},
\end{equation}
where $\hat{f}(\cdot)$ denotes the concave envelope of $f(\cdot)$ over the domain $\mathcal{X}$, i.e., $\hat{f} = -\operatorname{conv}(-f)$.
Therefore, the task of convexifying $\mathcal{Z}$ reduces entirely to characterizing $\hat{f}$. 
In Section~\ref{sec:application_qp}, we demonstrate how some mixed-integer quadratically constrained sets can be reformulated into $\mathcal{Z}$. 
Although obtaining a complete characterization of the concave envelope $\hat{f}$ is challenging, we show how a strong relaxation can be derived using the concave envelope of $f(\cdot)^2$ and establish an approximation gap. 
We assess the effectiveness of this reformulation technique through a computational study on a distributionally robust chance-constrained mixed-binary (multi-dimensional) knapsack problem. 
The results indicate that our reformulation yields a tighter bound and significantly reduces the solution time compared to a naive formulation.

\noindent\textbf{Notations.} Throughout this paper, we use boldface letters to denote vectors. Superscripts are used to distinguish between different vectors in a collection, while subscripts indicate specific components of a vector.
Given a matrix $A$, $\col(A)$ denotes the linear space spanned by the column vectors of $A$. Given a set $S$, $\proj_x(S)$ denotes its orthogonal projection onto the $x$-space. We let $\mathbb{N}$ denote the set of positive integers, i.e., $\mathbb{N}=\{1,2,\ldots\}$. For $n\in\mathbb{N}$, we let $\Delta_n$ denote the standard $n$-simplex in $\bbR^{n+1}$, i.e., $\Delta_n = \{\veclambda \in \bbR^{n+1}_+ : \sum_{k=1}^{n+1} \lambda_k = 1\}$.

\section{Ideal Formulation for $\calZ$}\label{sec:theory}

Without loss of generality, we may assume that $\calZ$ satisfies the following assumptions, possibly after a reformulation.
\begin{assumption}\label{assump:f}
    For all $\vecx\in\calX$, there exists $\vecy\in\bbR^{m}$, which may depend on $\vecx$, such that $\norm{A\vecx + B\vecy + \vecd}_2 \leq f(\vecx)$, i.e., $\proj_{\vecx}(\calZ)=\calX$.
\end{assumption}
\begin{assumption}\label{assump:col_space}
    All columns of $A$ and the vector $\vecd$ can be expressed as linear combinations of columns of $B$, i.e., $\col(A) \subseteq \col(B)$, and $\vecd\in \col(B)$.
\end{assumption}
Indeed, given $\calZ$, we may redefine $\calX$ as $\proj_{\vecx}(\calZ)$ so that Assumption \ref{assump:f} is satisfied. Regarding Assumption \ref{assump:col_space}, given $(A,B,\vecd)$ and $f(\cdot)$, one can decompose $A$ and $\vecd$ into $A = \bar{A} + A^\perp$ and $\vecd = \bar{\vecd} + \vecd^\perp$ such that $\col(\bar{A})\subseteq\col(B)$, $\bar{\vecd}\in\col(B)$, columns of $A^\perp$ and $\vecd^\perp$ lie in the kernel of $B$ by projecting the columns of $A$ and $\vecd$ onto $\col(B)$. Then we can rewrite $\calZ$ equivalently as the set of $(\vecx, \vecy)\in\calX\times\bbR^m$ such that 
\begin{equation}\label{eq:reform_soc}
    \norm{\bar{A}\vecx + B\vecy + \bar{\vecd}}_2 \leq f'(\vecx):=\sqrt{f^2(\vecx) - \norm{A^\perp \vecx + \vecd^\perp}_2^2},
\end{equation}
in which case Assumption \ref{assump:col_space} is satisfied by $(\bar{A},B,\bar{\vecd)}$ and $f'(\cdot)$. Note that $f':\calX\rightarrow\bbR$ is well-defined and nonnegative on $\calX$ under Assumption \ref{assump:f}.

Under Assumption 1, $f(\cdot)$ is nonnegative over its domain $\calX$, and so is $\hat{f}(\cdot)$ over $\conv(\calX)$. Our next result shows that one can characterize the convex hull of the set $\calZ$ given the concave envelope $\hat{f}(\cdot)$. 
\begin{theorem}\label{thm:conv_characterization}
    Suppose Assumptions \ref{assump:f} and \ref{assump:col_space} hold. Then the convex hull of $\calZ$ can be characterized as
    \begin{equation}\label{eq:cvx_characterization}
        \conv(\calZ) = \left\{(\vecx,\vecy) \in \conv(\calX) \times \bbR^m: \norm{A\vecx + B\vecy + \vecd}_2 \leq \hat{f}(\vecx) \right\},
    \end{equation}
    with $\conv(\calZ)$ being a closed set.
\end{theorem}

Before proving the theorem, the following lemma on the property of concave envelopes over compact sets is needed. The proof is presented in Appendix \ref{apdx:proof}. 
\begin{lemma}\label{lem:concave_envelope}
    Let $\calX \subsetneq \bbR^n$ be a compact set, $f: \calX \to \bbR$ be an upper semi-continuous function, and $\hat{f}:\conv(\calX)\rightarrow\bbR$ be its concave envelope over $\conv(\calX)$. Then, for any $\vecx \in \conv(\calX)$, there exist $K \leq n+1$ points $\vecx^1, \ldots, \vecx^K \in \calX$ and coefficients $\lambda_1, \ldots, \lambda_K > 0$ with $\sum_{k=1}^K \lambda_k = 1$ such that $\vecx = \sum_{k=1}^K \lambda_k \vecx^k$ and $\hat{f}(\vecx) = \sum_{k=1}^K \lambda_k f(\vecx^k)$. Moreover, $\hat{f}(\cdot)$ is upper semi-continuous on $\conv(\calX)$.
\end{lemma}

Now we are ready to prove Theorem \ref{thm:conv_characterization}.
\begin{proof}[of Theorem \ref{thm:conv_characterization}]
    Let $\calW$ be as defined in \eqref{eq:set_W}.
    Note that $\calW$ is closed and convex since $\hat{f}(\cdot)$ is upper semi-continuous and concave, and $\conv(\calX)\times \bbR^m$ is closed and convex. Therefore, we only need to show that $\conv(\calZ) = \calW$. We have $\calZ \subseteq \calW$ as by definition $\hat{f}(\vecx)\geq f(\vecx)$ for all $\vecx\in \calX$. Then the inclusion $\conv(\calZ)\subseteq \calW$ follows from convexity of $\calW$. It remains to show that $\calW\subseteq\conv(\calZ)$.
    
    We first consider the case when $B$ has full column rank, i.e., $B\vecy = \boldsymbol{0}$ if and only if $\vecy = \boldsymbol{0}$. Then $\calZ$ (resp. $\calW$) is closed and bounded and hence compact, as $f(\cdot)$ (resp. $\hat{f}(\cdot)$) is bounded and upper semi-continuous.
    By \cite{rockafellarConvexAnalysis1997}[Corollary 18.5.1], in order to prove $\calW \subseteq \conv(\calZ)$, it suffices to show that every extreme point of $\calW$ lies in $\calZ$. Let $(\bar{\vecx},\bar{\vecy})$ be an arbitrary extreme point of $\calW$. Then by Lemma \ref{lem:concave_envelope}, there exist $1\leq K \leq n+1$ points $\vecx^1, ..., \vecx^K \in \calX$ and $\lambda_1, ..., \lambda_K>0$ with $\sum_{k=1}^K \lambda_k = 1$, such that $\bar{\vecx} = \sum_{k=1}^K \lambda_k \vecx^k$ and $\hat{f}(\bar{\vecx}) = \sum_{k=1}^K \lambda_k f(\vecx^k)$. Then we have the following two cases.
    \begin{enumerate}
        \item [(a)] $K = 1$. In this case, we have $\bar{\vecx} = \vecx^1$ and $\hat{f}(\bar{\vecx}) = f(\vecx^1)$, implying $(\bar{\vecx}, \bar{\vecy}) \in \calZ$.
        \item [(b)] $K\geq 2$. We will prove by contradiction that $(\bar{\vecx}, \bar{\vecy})$ cannot be an extreme point of $\calW$. The key idea is finding two distinct points whose convex combination is $(\bar{\vecx}, \bar{\vecy})$, and the discussion is separated into the following three subcases. 
        \begin{enumerate}
            \item [(i)] $\norm{A\bar{\vecx} + B\bar{\vecy} +\vecd}_2 < \hat{f}(\bar{\vecx})$. For any nonzero vector $\vecu \in \bbR^m$, we can choose $\epsilon > 0$ small enough such that $\norm{A\bar{\vecx} + B(\bar{\vecy}\pm \epsilon \vecu) + \vecd}_2 < \hat{f}(\bar{\vecx})$. Thus $(\bar{\vecx}, \bar{\vecy}\pm \epsilon \vecu) \in \calW$, contradicting the fact that $(\bar{\vecx}, \bar{\vecy})$ is an extreme point of $\calW$.
            \item [(ii)] $\norm{A\bar{\vecx} + B\bar{\vecy} +\vecd}_2 = \hat{f}(\bar{\vecx}) = 0$. Then we must have $f(\vecx^k) = 0$ for $k = 1,...,K$, since $f(\vecx^k)\geq 0$ by Assumption \ref{assump:f} and $\sum_{k=1}^{K} \lambda_k f(\vecx^k) =\hat{f}(\bar{\vecx}) = 0$ in this case. Let $\vecv = \vecx^1 - \bar{\vecx}$. Then we can perturb $\bar{\vecx}$ along the direction and the reverse direction of $\vecv$ while keeping it in $\conv(\{\vecx^1, \ldots, \vecx^K\})\subseteq\conv(\calX)$ since $K\geq 2$. Let $\vecu \in \bbR^m$ satisfy $B\vecu = -A\vecv$. Such $\vecu$ exists since $\col(A)\subseteq \col(B)$. For $\epsilon > 0$ small enough (so that $\bar{\vecx} \pm \epsilon \vecv$ are still in $\conv\{\vecx^1, ..., \vecx^K\}$)
            \begin{align*}
                &\norm{A(\bar{\vecx}\pm \epsilon \vecv) + B(\bar{\vecy}\pm \epsilon\vecu) + \vecd}_2 \\ 
                =& \norm{A\bar{\vecx}+B\bar{\vecy}+\vecd}_2\\ 
                =& 0 \\ 
                \leq & \hat{f}(\bar{\vecx}\pm \epsilon \vecv).
            \end{align*}
            The last inequality follows from nonnegativity of $\hat{f}$ (due to Assumption \ref{assump:f}) and the fact that $\bar{\vecx}\pm \epsilon\vecv \in \conv(\{\vecx^1, \ldots, \vecx^K\})\subseteq\conv(\calX)$. Hence we have $(\bar{\vecx}\pm \epsilon\vecv, \bar{\vecy}\pm \epsilon\vecu) \in \calW$, contradicting to the fact that $(\bar{\vecx}, \bar{\vecy})$ is an extreme point.
            \item [(iii)] $\norm{A\bar{\vecx} + B\bar{\vecy} +\vecd}_2 = \hat{f}(\bar{\vecx}) > 0$. Let $\vecv$ be as defined above. But $\vecu$ in this case has to be chosen more carefully since the value of $\hat{f}(\cdot)$ may go down around $\bar{\vecx}$ in this case. Specifically, by Assumption \ref{assump:col_space}, we let $\vecu\in \bbR^m$ be chosen such that 
            \begin{equation*}
                B\vecu = \frac{f(\vecx^1) - \hat{f}(\bar{\vecx})}{\hat{f}(\bar{\vecx})}(A\bar{\vecx}+B\bar{\vecy} + \vecd) - A\vecv.
            \end{equation*}
            Then, for small enough $\epsilon>0$, we have 
            \begin{align}
                &\norm{A(\bar{\vecx}+\epsilon\vecv) + B(\bar{\vecy}+\epsilon \vecu) +\vecd}_2 \nonumber\\
                =& \norm{\frac{(1-\epsilon)\hat{f}(\bar{\vecx}) + \epsilon f(\vecx^1)}{\hat{f}(\bar{\vecx})}(A\bar{\vecx} + B\bar{\vecy} + \vecd)}_2\nonumber\\
                =& (1-\epsilon)\hat{f}(\bar{\vecx}) + \epsilon f(\vecx^1)\nonumber\\
                \leq& (1-\epsilon)\hat{f}(\bar{\vecx}) + \epsilon \hat{f}(\vecx^1) \leq \hat{f}(\bar{\vecx} + \epsilon\vecv),\label{ineq:plus}
            \end{align}
            where the second equality follows from positive homogeneity of the norm and the fact that we can choose $\epsilon$ small enough so that $(1-\epsilon)\hat{f}(\bar{\vecx}) + \epsilon f(\vecx^1)>0$, 
            and the last inequality follows from concavity of $\hat{f}(\cdot)$. On the other hand, for small enough $\epsilon>0$, we have
            \begin{align}
                &\norm{A(\bar{\vecx}-\epsilon\vecv) + B(\bar{\vecy}-\epsilon \vecu) +\vecd}_2 \nonumber\\
                =& \norm{\frac{(1+\epsilon)\hat{f}(\bar{\vecx}) - \epsilon f(\vecx^1)}{\hat{f}(\bar{\vecx})}(A\bar{\vecx} + B\bar{\vecy} + \vecd)}_2\nonumber\\ 
                =& (1+\epsilon)\hat{f}(\bar{\vecx}) - \epsilon f(\vecx^1)\nonumber\\
                =& \left(1+\epsilon - \frac{\epsilon}{\lambda_1}\right)\hat{f}(\bar{\vecx}) + \frac{\epsilon}{\lambda_1}\sum_{k=2}^K \lambda_k f(\vecx^k) \nonumber\\ 
                \leq& \left(1+\epsilon - \frac{\epsilon}{\lambda_1}\right)\hat{f}(\bar{\vecx}) + \frac{\epsilon}{\lambda_1}\sum_{k=2}^K \lambda_k \hat{f}(\vecx^k) \nonumber\\
                \leq& \hat{f}\left(\left(1+\epsilon-\frac{\epsilon}{\lambda_1}\right)\bar{\vecx} + \frac{\epsilon}{\lambda_1}\sum_{k=2}^K\lambda_k \vecx^k\right) \nonumber\\ 
                =& \hat{f}(\bar{\vecx} - \epsilon\vecv),\label{ineq:minus}
            \end{align}
            where the second equality holds if we choose $\epsilon>0$ small enough so that $(1+\epsilon)\hat{f}(\bar{\vecx}) - \epsilon f(\vecx^1)>0$, and the third and the last equalities follow from the facts that $\hat{f}(\bar{\vecx}) = \sum_{k=1}^K \lambda_k f(\vecx^k)$ and $\bar{\vecx} = \sum_{k=1}^K \lambda_k \vecx^k$. 
            To conclude, by inequalities \eqref{ineq:plus} and \eqref{ineq:minus}, we have $(\bar{\vecx}\pm\epsilon\vecv, \bar{\vecy}\pm \epsilon\vecu)\in \calW$ for small enough $\epsilon>0$, whose midpoint is $(\bar{\vecx}, \bar{\vecy})$.
        \end{enumerate}
    \end{enumerate}
    In summary, we have demonstrated that any extreme point of $\calW$ must belong to $\calZ$, and thus $\calW = \conv(\calZ)$, when $B$ has full column rank, i.e., $\rank(B)=m$.
    
    Now consider the case when $\rank(B) = \hat{m} <m$. In this case, we can find a full-rank matrix $U\in \bbR^{m\times m}$ (corresponding to some elementary column operation) such that $BU = (\hat{B}, \boldsymbol{0})$, where $\hat{B}$ has $\hat{m}$ columns with full column rank. Define 
    \begin{align*}
        \hat{\calW} &= \{(\vecx, \vecw) \in \conv(\calX)\times \bbR^m: \norm{A\vecx + \hat{B}\vecw_{1:\hat{m}} + \vecd}_2\leq \hat{f}(\vecx)\}, \\ 
        \hat{\calZ} &=  \left\{(\vecx,\vecw) \in \calX \times \mathbb{R}^m: \norm{A\vecx + \hat{B}\vecw_{1:\hat{m}} + \vecd}_2 \leq f(\vecx) \right\},
    \end{align*}
    where $\vecw_{1:\hat{m}}$ is the subvector of $\vecw$ consisting of the first $\hat{m}$ elements. The two sets are unconstrained in the last $m - \hat{m}$ dimensions. Let $\Tilde{\calW}$ (resp. $\Tilde{\calZ}$) denote the projection of $\hat{\calW}$ (resp. $\hat{\calZ}$) onto the first $(n+\hat{m})$-dimensional subspace. Then we have the following observations: 
    \begin{enumerate}
        \item [(a)] $\left(\begin{array}{cc}
        I_n & \boldsymbol{0} \\
        \boldsymbol{0} & U^{-1}
        \end{array}\right):\bbR^{n+m} \to \bbR^{n+m}$ is a linear bijection from $\calW$ (resp. $\calZ$) to $\hat{\calW}$ (resp. $\hat{\calZ}$);
        \item [(b)] $\hat{\calW} = \Tilde{\calW}\times\bbR^{m- \hat{m}}$ and $\hat{\calZ} = \Tilde{\calZ}\times\bbR^{m- \hat{m}}$;
        \item [(c)] $\hat{\calW}$ is convex and $\conv(\hat{\calZ}) = \conv(\Tilde{\calZ}) \times \bbR^{m-\hat{m}}$;
        \item [(d)] $\Tilde{\calW}$ and $\Tilde{\calZ}$ are compact sets respectively admitting characterization in the form of \eqref{eq:set_W} and \eqref{eq:set_Z} with $\hat{B}$ having full column rank and Assumptions \ref{assump:f} and \ref{assump:col_space} holding.
    \end{enumerate}
    Applying our proof for the full-rank $B$ case, we know that $\Tilde{\calW} = \conv(\Tilde{\calZ})$ and they are compact. Therefore, $\hat{\calW} = \conv(\hat{\calZ})$ and they are closed. It then follows from observation (a) that $\calW = \conv(\calZ)$ and they are closed, as linear bijection preserves closedness and convexity.\qed
\end{proof}

Note that the above proof relies only on positive homogeneity and nonnegativity of the 2-norm, properties that hold for any norm. Therefore, Theorem \ref{thm:conv_characterization} remains valid if the 2-norm in \eqref{eq:set_Z} and \eqref{eq:cvx_characterization} is replaced by any other norm. 
\begin{corollary}
    Let $\norm{\cdot}$ be any norm defined in the Euclidean space. Define $\calZ$ as
    $$
    \calZ := \left\{(\vecx,\vecy) \in \calX \times \mathbb{R}^m: \norm{A\vecx + B\vecy + \vecd} \leq f(\vecx) \right\}.
    $$
    Suppose Assumptions \ref{assump:f} and \ref{assump:col_space} hold.  Then the convex hull of $\calZ$ can be characterized as
    $$
        \conv(\calZ) = \left\{(\vecx,\vecy) \in \conv(\calX) \times \bbR^m: \norm{A\vecx + B\vecy + \vecd} \leq \hat{f}(\vecx) \right\},
    $$
    with $\conv(\calZ)$ being a closed set.
\end{corollary}
However, we shall also note that if Assumption 2 is not originally satisfied by $\calZ$, then in general there may not exist a reformulation of $\calZ$ like \eqref{eq:reform_soc} in the 2-norm case such that Assumption 2 holds after reformulation.

Assumption \ref{assump:col_space} ensures that, whichever direction we perturb $\vecx$ there exists a perturbation of $\vecy$ such that $A\vecx+B\vecy + \vecd$ is scaled with an arbitrary positive factor we choose, while keeping the direction unchanged. That depicts the essential property of the conic set $\calZ$ that underlies Theorem \ref{thm:conv_characterization}. However, if $\vecy$ is constrained, this property may no longer hold due to the reduced degrees of freedom in $\vecy$, resulting in the characterization \eqref{eq:cvx_characterization} inexact. The following two examples illustrate the scenarios where Assumption \ref{assump:col_space} is violated and $\vecy$ is constrained, respectively.

\begin{figure}[tp]
\begin{subfigure}{.45\textwidth}
\centering
\includegraphics[width=\textwidth]{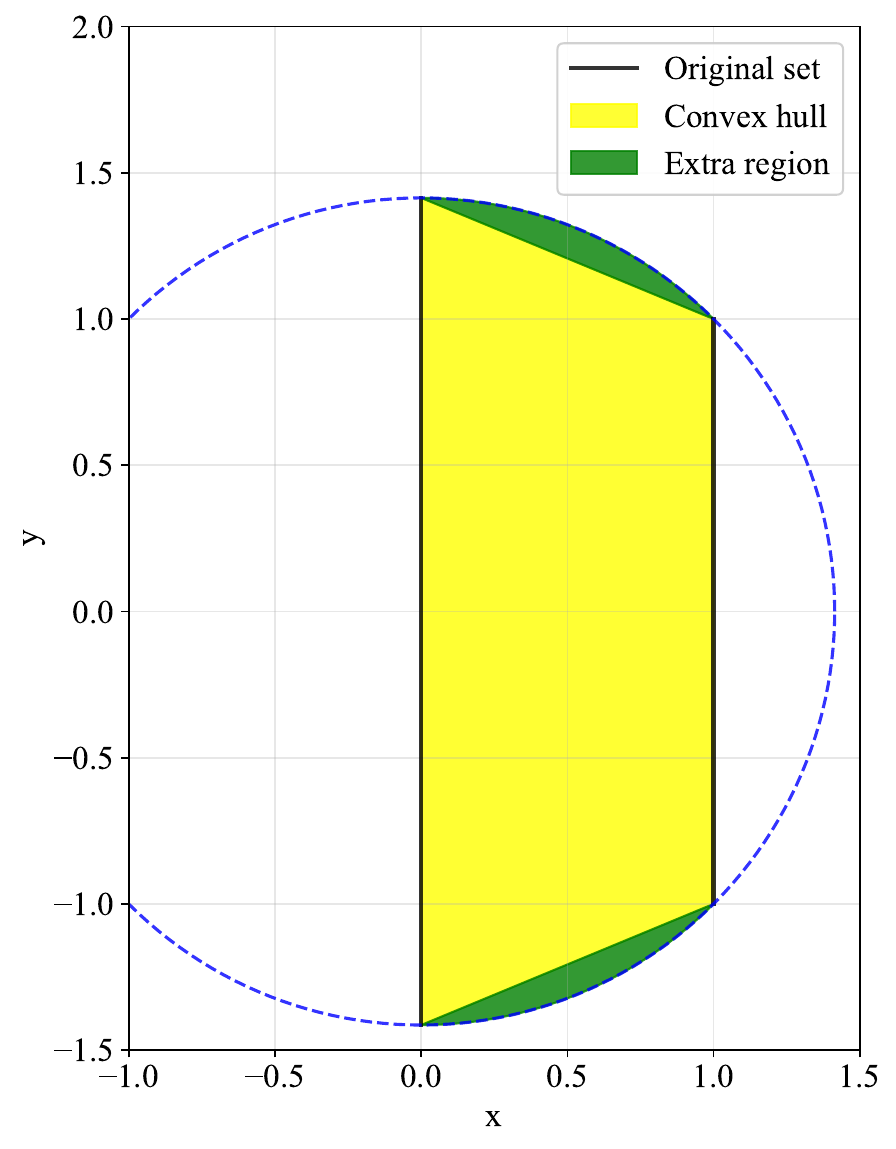}
\caption{Example \ref{eg:assumption2}}
\label{fig:assumption2}
\end{subfigure}
\hfill
\begin{subfigure}{.45\textwidth}
\centering
\includegraphics[width=\textwidth]{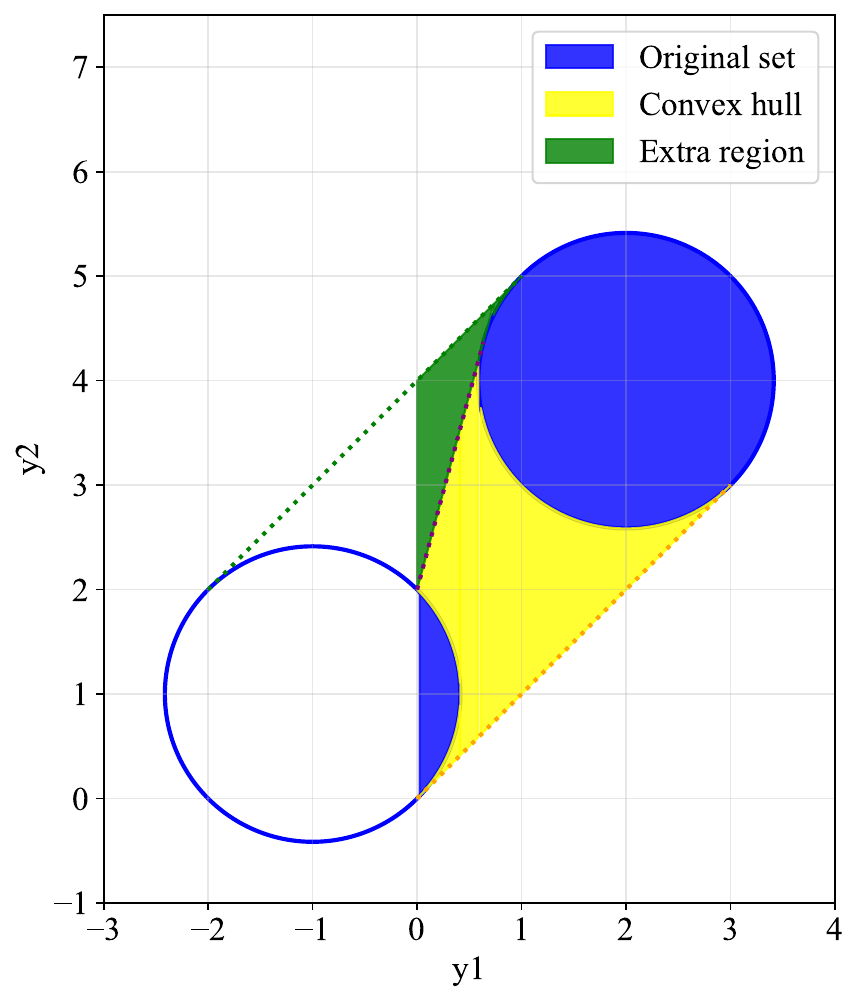}
\caption{Example \ref{eg:y_constrained}}
\label{fig:y_constrained}
\end{subfigure}
\caption{The convex hull of $\calZ$ (a) when Assumption \ref{assump:col_space} does not hold, or (b) when 
$\vecy$ is further constrained.} \label{fig:egs}
\end{figure}
\begin{example} \label{eg:assumption2}
    Let $m = n = 1$, $\calX = \{0,1\}$, $A = (1,0)\trans$, $B = (0,1)\trans$, $\vecd = \boldsymbol{0}$, and $f(x) = \sqrt{2}$ be a constant function over $\calX$. Then set $\calZ$ is given by
    $$\calZ = \left\{(x,y)\in \{0,1\}\times \bbR : \sqrt{x^2 + y^2} \leq \sqrt{2}\right\}.$$
    For this particular $f(\cdot)$, we have $\hat{f}(x) = \sqrt{2}$ for all $x\in\conv(\calX) = [0,1]$ while Assumption \ref{assump:col_space} is not satisfied. In this case, $\calW$ actually includes some extra region apart from $\conv(\calZ)$, as shown in  Figure \ref{fig:assumption2}. To obtain the $\conv(\calZ)$, following \eqref{eq:reform_soc}, one may first reformulate $\calZ$ as\begin{align*}
        \calZ = \left\{(x,y)\in \{0,1\}\times \bbR : |y| \leq f'(x):=\sqrt{2-x^2}\right\}
    \end{align*}
    so that Assumption \ref{assump:col_space} is satisfied,
    and note that $\hat{f}'(x)=(1-\sqrt{2})x+\sqrt{2}$ over $\conv(\calX)=[0,1]$. Then by Theorem \ref{thm:conv_characterization} we have
    \begin{align*}
        \conv(\calZ) = \left\{(x,y)\in [0,1]\times \bbR : |y| \leq (1-\sqrt{2})x+\sqrt{2}\right\}.
    \end{align*}
    \hfill$\blacksquare$
\end{example}

\begin{example} \label{eg:y_constrained}
    Let $m = 1$ and $n = 2$, and define $\calZ$ as 
    \begin{equation}\label{exp:constrained}
        \calZ = \left\{(x,y_1,y_2) \in \{0,1\} \times \bbR^2_+ : \norm{\left(\begin{array}{c}
            3 \\
            3
        \end{array}\right)x + \left(\begin{array}{c}
            y_1 \\
            y_2
        \end{array}\right) - \left(\begin{array}{c}
            -1 \\
            1
        \end{array}\right)}_2 \leq \sqrt{2}\right\}.
    \end{equation}
    Note that, in this example, nonnegativity constraints are imposed on continuous variables $(y_1,y_2)$. As shown in Figure \ref{fig:y_constrained}, the nonnegative constraint cuts off part of the left ball, and $\conv(\calZ)$ is strictly smaller than the set obtained by directly convexifying $\calX$ in \eqref{exp:constrained}.
    \hfill$\blacksquare$
\end{example}

\section{Applications to Mixed-Integer Quadratically Constrained Programming}\label{sec:application_qp}
One of the primary reasons for investigating $\calZ$ is its frequent occurrence in Mixed-Integer Quadratically Constrained Programs (MIQCPs) as a substructure. Indeed, by properly choosing $(A,B,\vecd)$ and $f(\cdot)$, any constraint of the form\begin{equation}\label{ineq:quad}
    \left(\begin{array}{c}
        \vecx \\
        \vecy
    \end{array}\right) \trans
    \left(\begin{array}{cc}
        Q_{xx} & Q_{xy} \\
        Q_{xy}\trans & Q_{yy}
    \end{array}\right)
    \left(\begin{array}{c}
        \vecx \\
        \vecy
    \end{array}\right) + 
    \left(\begin{array}{c}
        \veca_x \\
        \veca_y
    \end{array}\right)\trans
    \left(\begin{array}{c}
        \vecx \\
        \vecy
    \end{array}\right) \leq g(\vecx)
\end{equation}
with a positive definite $Q_{yy}$ can be reformulated as $\norm{A\vecx + B\vecy + \vecd}_2 \leq f(\vecx)$, as demonstrated in Appendix \ref{apdx:reform_qc}. In the most common case when $g(\cdot)$ is a constant, $f(\cdot)$ can be expressed as $\sqrt{q(\cdot)}$ for some quadratic function $q(\cdot)$. Then, by applying Theorem \ref{thm:conv_characterization}, one can reduce the convexification of the substructure \eqref{ineq:quad} (together with some constraints on $\vecx$) to the convexification of the hypograph of $f(\cdot)=\sqrt{q(\cdot)}$ (restricted to some domain). 
While the convexification of the epi-/hypo-graph of a quadratic function over integer or binary points has been widely studied in the literature, notably on the Boolean quadric polytope (BQP) \cite{boros1993cut,de1990cut,padberg1989boolean}, the convexification of the square root of a quadratic function over such points remains largely unexplored (with an exception in the submodular case \cite{atamturk2020submodularity,atamturk2008polymatroids,yu2017polyhedral,zhangAmbiguousChanceConstrainedBinary2018}). However, as shown in the next result, one can utilize the knowledge about the concave envelope $\hat{q}(\cdot)$ of a function $q(\cdot)$ to derive convex outer approximations for the concave envelope of $f(\cdot)=\sqrt{q(\cdot)}$. The proof appears in Appendix \ref{apdx:proof}. 
\begin{proposition}\label{prop:relax_gap}
    Let $\calX\subsetneq \bbR^n$ be a compact set, $f: \calX \to \bbR_+$ be an upper semi-continuous function, and $\hat{f}:\conv(\calX)\rightarrow\bbR_+$ be its concave envelope over $\conv(\calX)$. Let $q:\calX\rightarrow\bbR_+$ be such that $q(\vecx)=f(\vecx)^2$ for all $\vecx\in\calX$, and $\hat{q}:\conv(\calX)\rightarrow\bbR_+$ be its concave envelope over $\conv(\calX)$. Suppose $f(\vecx)\in[L,U]$ for all $\vecx\in\calX$. Then $\sqrt{\hat{q}(\cdot)}$ is concave and upper semi-continuous over $\conv(\calX)$, and for all $\vecx\in\conv(\calX)$, we have \begin{align}\label{ineq:bound}
        0\leq\sqrt{\hat{q}(\vecx)}-\hat{f}(\vecx)\leq\frac{(U-L)^2}{4(L+U)},
    \end{align}
    where $\frac{0}{0}:=0$.
\end{proposition}

Proposition \ref{prop:relax_gap} motivates the following convex relaxation of $\conv(\calZ)$:\begin{align}
    \calR(\calZ):=\left\{(\vecx,\vecy) \in \conv(\calX) \times \bbR^m: \norm{A\vecx + B\vecy + \vecd}_2 \leq \sqrt{\hat{q}(\vecx)} \right\}.
\end{align}
In practice, in order to obtain $\calR(\calZ)$ or approximations thereof, we propose the following reformulation of $\calZ$ in the extended variable space,\begin{align}
   \calZ=\proj_{(\vecx,\vecy)}\left\{(\vecx,\vecy,\eta,\tau) \in \calX\times \bbR^{m}\times\bbR^2_+:\right.\nonumber& \\
    \norm{A\vecx + B\vecy + \vecd}_2 \leq \eta,&~\left.\eta^2\leq\tau,~\tau\leq {q}(\vecx)\right\},\label{constr:ideal_soc_qp}
\end{align}
to exploit standard MIP solvers' ability to approximate the concave envelope of $q(\cdot)$ using cutting planes in possibly an even further extended variable space, especially in the case when $q(\cdot)$ is quadratic. Note that Assumption \ref{assump:f} may not initially hold, in which case we may also rely on solver cutting planes to characterize $\conv(\proj_{\vecx}(\calZ))$.

\section{Computational Study}\label{sec:numerical_study}
In this section, we evaluate the effectiveness of the reformulation proposed in Section \ref{sec:application_qp} on a distributionally robust chance-constrained mixed-binary (multi-dimensional) knapsack problem (DRCC-MB-(M)KP), 
given by
\begin{subequations}\label{DRCC-KP}
    \begin{align}
        \max \textrm{ } & \vecp_x\trans \vecx + \vecp_y\trans \vecy \\
        \textrm{s.t. } & \inf_{\bbP \in \calD_j}\left\{\bbP((\vecw_x^j)\trans \vecx + (\vecw_y^j)\trans \vecy \leq c^j)\right\}\geq 1-\alpha, \quad j\in \calJ \label{constr:cc}\\
        & \vecx\in \{0,1\}^n, \vecy\in [0,1]^m,
    \end{align}
\end{subequations}
where $n$ denotes the number of discrete items, $m$ the number of continuous items, and $\calJ$ indexes the multiple resource types. The resource consumptions and profits of the discrete items are denoted by vectors $\vecw_x\in\bbR^n$ and $\vecp_x\in\bbR^n$, while those of the continuous items are denoted by $\vecw_y\in\bbR^m$ and $\vecp_y\in\bbR^m$, respectively. The capacity limit of each resource type $j\in\calJ$ is denoted by $c^j$.
For each resource type $j\in\calJ$, the distributionally robust chance constraint (DRCC) \eqref{constr:cc} enforces that the total consumption shall not exceed its capacity limit with probability at least $1-\alpha$ under any distribution of resource consumption $(\vecw_x^j,\vecw_y^j)$ in the ambiguity set $\calD_j$. The ambiguity set $\calD_j$ is chosen to be a Chebyshev ambiguity set, i.e., the set of all distributions with some given mean and covariance matrix. Then, for each resource type $j\in \calJ$, it is known \cite{ghaoui2003worst} that the DRCC \eqref{constr:cc} is equivalent to
\begin{equation} \label{constraint:soc}
     \left(\begin{array}{c}\vecmu_x\\ \vecmu_y \end{array} \right)\trans\left(\begin{array}{c}\vecx\\ \vecy \end{array} \right)
    + \sqrt{\frac{1-\alpha}{\alpha}} 
    \sqrt{\left(\begin{array}{c}\vecx\\ \vecy \end{array}\right)\trans \left(\begin{array}{cc}
        \Sigma_{xx} & \Sigma_{xy} \\
        \Sigma_{xy}\trans & \Sigma_{yy}
    \end{array}\right) \left(\begin{array}{c}\vecx\\ \vecy \end{array}\right) } \leq c,
\end{equation}
where the dependence on the resource type $j$ is omitted for brevity, and $\vecmu$ and $\Sigma$ denote the given mean and covariance of resource consumptions, respectively.
We may further reformulate constraint \eqref{constraint:soc} as 
\begin{subequations}\label{constr:soc}
    \begin{align}
        \left(\begin{array}{c}\vecx\\ \vecy \end{array}\right)\trans  \left(\begin{array}{cc}
        \Tilde{\Sigma}_{xx} & \Tilde{\Sigma}_{xy} \\
        \Tilde{\Sigma}_{xy}\trans & \Tilde{\Sigma}_{yy}
    \end{array}\right)\left(\begin{array}{c}\vecx\\ \vecy \end{array}\right) +& 2c\left(\begin{array}{c}\vecmu_x\\ \vecmu_y \end{array} \right)\trans\left(\begin{array}{c}\vecx\\ \vecy \end{array} \right) \leq c^2, \label{ineq:cc_quad} \\
    \vecmu_x\trans \vecx + \vecmu_y\trans \vecy\leq& c,\label{ineq:linear}
    \end{align}
\end{subequations}
where $\Tilde{\Sigma}_{xx} = \Tilde{\alpha}\Sigma_{xx} - \vecmu_x\vecmu_x\trans$, $\Tilde{\Sigma}_{xy} = \Tilde{\alpha}\Sigma_{xy} - \vecmu_x \vecmu_y\trans$, $\Tilde{\Sigma}_{yy} = \Tilde{\alpha}\Sigma_{yy} - \vecmu_y\vecmu_y\trans$, and $\Tilde{\alpha} = (1-\alpha)/\alpha$. Once $\Sigma_{yy} \succ 0$ and $\alpha$ is small enough, we have $\Tilde{\Sigma}_{yy} \succ 0$. As discussed in Section \ref{sec:application_qp}, each constraint of \eqref{ineq:cc_quad} can be reformulated as $\norm{A\vecx + B\vecy + \vecd}_2 \leq f(\vecx)$, and hence into the form of \eqref{constr:ideal_soc_qp}. In the remainder of this section, we test the following two MIQCP formulations for DRCC-MB-(M)KP:\begin{enumerate}
    \item {\CCP}: Model \eqref{DRCC-KP} with all DRCCs \eqref{constr:cc} reformulated into the form of \eqref{constr:soc};
    \item {\SOC}: Model \eqref{DRCC-KP} with all DRCCs \eqref{constr:cc} reformulated into the form of \eqref{constr:ideal_soc_qp} (using auxiliary variables) and \eqref{ineq:linear}.
\end{enumerate}
Formulation {\CCP} represents the vanilla MIQCP formulation of the original problem \eqref{DRCC-KP}.
Formulation {\SOC} explicitly exposes the SOC structure and the BQP structure, which enables the solver to generate cutting planes to approximate the hypographs of quadratic functions. It is observed that the explicit hypograph structure $\eta\leq q(\vecx)$ is essential to trigger Gurobi's ability to generate such cuts.

We modify existing instances of deterministic (multi-dimensional) knapsack problems \cite{martelloDynamicProgrammingStrong1999,petersenComputationalExperienceVariants1967} to construct our DRCC-MB-(M)KP instances.
The construction of DRCC-MB-(M)KP instances and experiment settings are detailed in Appendix \ref{apdx:instance_generation} and \ref{apdx:exp_settings}, respectively.


\begin{table}[tb!]
    \centering
    \small
    \caption{The average root gap, number of solved instances, and the average optimality gap over different types of DRCC-MB-KP problems with varying problem sizes.}
    \label{tab:gap_kp}
\begin{threeparttable}
\begin{tabular}{ccrrrrrr}
    \toprule
\multirow{2}{*}{$(n,m)$} & \multirow{2}{*}{Type} & \multicolumn{2}{c}{Root Gap} & \multicolumn{2}{c}{\# Solved} & \multicolumn{2}{c}{Optimality Gap}   \\
\cmidrule(lr){3-4} \cmidrule(lr){5-6} \cmidrule(lr){7-8}
 &   & \CCP & \SOC  & \CCP & \SOC & \CCP & \SOC   \\
\midrule \midrule \addlinespace
\multirow{4}{*}{(\phantom{0}25, \phantom{0}25)}  & 1 & 695.9\%  & \textbf{272.5\%}  & 4/5& \textbf{5/5}     & 409.7\%         & \textbf{0.0\%}   \\
      & 2 & 450.2\%  & \textbf{113.7\%}  & 3/5& \textbf{5/5}     & 22.7\%          & \textbf{0.0\%}   \\
      & 3 & 431.0\%  & \textbf{66.3\%}   & \textbf{5/5}     & \textbf{5/5}     & \textbf{0.0\%}  & \textbf{0.0\%}   \\
      & 4 & 303.0\%  & \textbf{43.2\%}   & \textbf{5/5}     & \textbf{5/5}     & \textbf{0.0\%}  & \textbf{0.0\%}   \\
\cmidrule(lr){1-8}
\multirow{4}{*}{(\phantom{0}50, \phantom{0}50)}  & 1 & 250.8\%  & \textbf{134.4\%}  & 0/5& \textbf{4/5}     & 309.2\%         & \textbf{2.2\%}   \\
      & 2 & 427.4\%  & \textbf{273.3\%}  & 0/5& \textbf{2/5}     & 411.6\%         & \textbf{3.3\%}   \\
      & 3 & 394.4\%  & \textbf{235.1\%}  & 0/5& \textbf{2/5}     & 422.9\%         & \textbf{3.6\%}   \\
      & 4 & 260.8\%  & \textbf{109.4\%}  & 2/5& \textbf{4/5}     & 14.2\%          & \textbf{1.0\%}   \\
\cmidrule(lr){1-8}
\multirow{4}{*}{(100, 100)} & 1 & 150.6\%  & \textbf{86.7\%}   & (1) 0/5           & 0/5& 160.2\%         & \textbf{63.5\%}  \\
      & 2 & 284.6\%  & \textbf{210.6\%}  & 0/5& 0/5& 333.3\%         & \textbf{189.2\%} \\
      & 3 & 294.8\%  & \textbf{210.9\%}  & (3) 0/5           & 0/5& 395.9\%         & \textbf{198.3\%} \\
      & 4 & 229.4\%  & \textbf{175.6\%}  & 0/5& 0/5& 263.1\%         & \textbf{154.1\%}\\
\bottomrule
\end{tabular}
\begin{tablenotes}
    \scriptsize
    \item The numbers in parentheses indicate the number of instances that terminate prematurely due to numerical issues.
\end{tablenotes}
\end{threeparttable}
\end{table}
\begin{table}[tb!]
    \centering
    \small
    \caption{The average times and the average numbers of explored nodes over different types of DRCC-MB-KP problems with varying problem sizes.}
    \label{tab:time_kp}
\begin{threeparttable}
\begin{tabular}{ccrrrrrr}
    \toprule
\multirow{2}{*}{$(n,m)$} & \multirow{2}{*}{Type} & \multicolumn{2}{c}{Time (s)} & \multicolumn{2}{c}{\# Nodes Explored}   \\
\cmidrule(lr){3-4} \cmidrule(lr){5-6}
 & & \CCP & \SOC & \CCP & \SOC \\
 \midrule \midrule \addlinespace
\multirow{4}{*}{(\phantom{0}25, \phantom{0}25)}  & 1 & 787.6 & \textbf{1.4}    & $ 9.7 \times 10^{4}$  & $ 6.6 \times 10^{2}$  \\
 & 2 & 1479.4     & \textbf{1.4}    & $ 1.2 \times 10^{5}$  & $ 8.8 \times 10^{2}$  \\
 & 3 & 423.3 & \textbf{5.3}    & $ 8.5 \times 10^{4}$  & $ 2.2 \times 10^{3}$  \\
 & 4 & 744.2 & \textbf{1.3}    & $ 7.9 \times 10^{4}$  & $ 8.3 \times 10^{2}$  \\
 \cmidrule(lr){1-6}
\multirow{4}{*}{(\phantom{0}50, \phantom{0}50)}  & 1 & 3600.2     & \textbf{931.5}  & $ 3.4 \times 10^{5}$  & $ 1.1 \times 10^{4}$  \\
 & 2 & 3600.1     & \textbf{2562.7} & $ 4.0 \times 10^{5}$  & $ 6.2 \times 10^{4}$  \\
 & 3 & 3600.1     & \textbf{2676.3} & $ 4.2 \times 10^{5}$  & $ 5.7 \times 10^{4}$  \\
 & 4 & 2196.8     & \textbf{1172.3} & $ 2.5 \times 10^{5}$  & $ 2.9 \times 10^{4}$  \\
  \cmidrule(lr){1-6}
\multirow{4}{*}{(100,100)} & 1 & (1) 2882.1  & 3600.1 & $ 1.8 \times 10^{5}$  & $ 7.4 \times 10^{2}$  \\
 & 2 & 3601.7     & 3600.1 & $ 2.5 \times 10^{5}$  & $ 6.1 \times 10^{2}$  \\
 & 3 & (3) 1441.0  & 3600.4 & $ 9.1 \times 10^{4}$  & $ 5.9 \times 10^{2}$  \\
 & 4 & 3600.6     & 3600.1 & $ 3.5 \times 10^{5}$  & $ 5.9 \times 10^{2}$
  \\
 \bottomrule
\end{tabular}
\begin{tablenotes}
    \scriptsize
    \item The numbers in parentheses indicate the number of instances that terminate prematurely due to numerical issues.
\end{tablenotes}
\end{threeparttable}    
\end{table}
We present the comparison of {\CCP} and {\SOC} on DRCC-MB-KP in Tables \ref{tab:gap_kp} and \ref{tab:time_kp}. For each type and size, we present the root gap, the number of instances solved to optimality, the ending optimality gap, the solution time, and the number of nodes explored before termination. The root gap is defined as the relative gap between the dual bound at the end of the root node (after cut generation) and the best objective value obtained by the two models within the time limit. Numbers in parentheses indicate the number of instances that terminate prematurely due to numerical issues. Each row is averaged over 5 instances.

As shown in Tables \ref{tab:gap_kp} and \ref{tab:time_kp}, {\SOC} consistently yields tighter root bounds and smaller optimality gaps compared to {\CCP} across all problem sizes and types. For small instances (i.e., with $50$ items), {\SOC} solves all instances within seconds, while {\CCP} solves $17$ out of $20$ instances within an hour. On medium instances (i.e., with $100$ items), {\SOC} achieves optimality on more than half of the instances, whereas {\CCP} generally fails to reach optimality, resulting in optimality gaps orders of magnitude larger than those of {\SOC}. 
{\SOC} demonstrates the ability to significantly tighten the optimality gap while exploring far fewer nodes. This efficiency stems from the automatic model lifting by linearization of the quadratic terms during the presolving phase and the generation of effective Gurobi cuts, primarily BQP cuts \cite{padberg1989boolean} and RLT cuts \cite{sheraliReformulationLinearizationTechniqueSolving1999}, which are absent in {\CCP}. However, for large instances (i.e., with 200 items), the abundance of these cuts and auxiliary linearization variables can significantly slow down the branch-and-bound process (orders of magnitude branching nodes explored within the time limit), leading to much less competitive performance of {\SOC}. 
Apart from solution efficiency, we observe that the reformulated model {\SOC} exhibits better numerical stability, as {\CCP} encountered four numerical issues while {\SOC} encountered none.

\begin{table}[tbp]
    \centering
    \small
    \caption{Root gaps, ending gaps, computation times, and number of branch-and-bound nodes over the seven DRCC-MB-MKP instances.}
    \label{tab:res_mkp}
\begin{threeparttable}
\begin{tabular}{crrrrrrrr}
    \toprule

\multirow{2}{*}{Instance} & \multicolumn{2}{c}{Root Gap}  & \multicolumn{2}{c}{Optimality Gap} & \multicolumn{2}{c}{Time (s)}  & \multicolumn{2}{c}{\# Nodes Explored}     \\
\cmidrule(lr){2-3} \cmidrule(lr){4-5} \cmidrule(lr){6-7} \cmidrule(lr){8-9}
   & \CCP & \SOC     & \CCP    & \SOC      & \CCP & \SOC   & \CCP & \SOC   \\
\midrule \midrule \addlinespace
1 & 688.3\% & \textbf{0.0\%}  & \textbf{*0.0\%} & \textbf{*0.0\%} & 0.0    & 0.0    & $ 1.7 \times10^{1}$ & $ 1.0 \times10^{0}$ \\
2 & \textbf{30.8\%} & 36.3\%  & \textbf{*0.0\%} & \textbf{*0.0\%} & 0.1    & 0.3    & $ 1.2 \times10^{2}$ & $ 1.2 \times10^{2}$ \\
3 & 20.5\%  & \textbf{11.8\%} & \textbf{*0.0\%} & \textbf{*0.0\%} & 0.3    & 0.5    & $ 7.6 \times10^{2}$ & $ 2.1 \times10^{2}$ \\
4 & 29.6\%  & \textbf{17.2\%} & \textbf{*0.0\%} & \textbf{*0.0\%} & 166.6  & \textbf{1.5}   & $ 4.0 \times10^{4}$ & $ 4.9 \times10^{2}$ \\
5 & 3.9\%   & \textbf{2.0\%}  & \textbf{*0.0\%} & \textbf{*0.0\%} & 31.4   & \textbf{4.4}   & $ 1.9 \times10^{4}$ & $ 1.2 \times10^{3}$ \\
6 & 99.1\%  & \textbf{20.9\%} & 140.2\% & \textbf{*0.0\%} & 3600.2 & \textbf{967.1} & $ 2.7 \times10^{5}$ & $ 7.0 \times10^{4}$ \\
7 & 70.2\%  & \textbf{19.8\%} & 103.9\% & \textbf{9.9\%}  & 3600.0 & 3600.0 & $ 1.8 \times10^{5}$ & $ 9.6 \times10^{4}$ \\
\bottomrule
\end{tabular}
\begin{tablenotes}
    \scriptsize
    \item [\textbf{*}] Indicates that the instance is solved to optimality.
\end{tablenotes}
\end{threeparttable}
\end{table}

Table \ref{tab:res_mkp} presents the computational results on seven DRCC-MB-MKP instances, where similar statistics are reported as Tables \ref{tab:gap_kp} and \ref{tab:time_kp}. The first three small-sized instances are solved to optimality by both methods within one second. For the relatively difficult instances (4 to 7), {\SOC} demonstrates much better performance. It yields tighter root bounds and solves instances 4 to 6 to optimality with shorter time and fewer branching nodes needed to reach optimality. In instance 7, {\SOC} terminates with an optimality gap $9.9\%$, significantly smaller than the gap $103.9\%$ for {\CCP}. 

\section{Conclusions}
In this paper, we study a class of mixed-integer conic sets and established that, under mild assumptions, their convex hull can be fully characterized via the concave envelope of the nonlinear function appearing on the right-hand side of the conic constraint.
Although explicitly deriving the convex hull is, in general, computationally challenging, the associated ideal formulation provides a foundation for strong relaxations through cutting-plane techniques. The effectiveness of these relaxations has been empirically validated in our numerical study.
Building on these findings, future research may focus on developing efficient procedures for generating cutting planes that approximate the convex envelope of square-root functions, particularly in the context of MIQCP, thereby enhancing solver performance on a broader class of mixed-integer nonlinear programs.
\bibliographystyle{splncs04}
\bibliography{ref}
\appendix
\section{Omitted Proofs}\label{apdx:proof}
\begin{proof}[of Lemma \ref{lem:concave_envelope}]
    By \cite{rockafellarConvexAnalysis1997}[Corollary 17.1.5], we have 
    $$
        \hat{f}(\vecx) = \sup\left\{ \sum_{k=1}^{n+1} \lambda_k f(\vecx^k) : \veclambda \in \Delta_n, (\vecx^k)_{k=1}^{n+1}\in \calX^{n+1}, \sum_{k=1}^{n+1} \lambda_k \vecx^k = \vecx\right\}.
    $$
    Note that the supremum is taken over a compact region, and the objective function is continuous in $(\veclambda, \vecx^1, \ldots, \vecx^{n+1})$. The supremum is thus attained at a point $(\veclambda, \vecx^1, ..., \vecx^{n+1})\in \Delta_n \times \calX^{n+1}$, which satisfies $\vecx = \sum_{k=1}^{n+1} \lambda_k \vecx^k$. The first half of the statement can then be obtained by dropping points $\vecx^k$ associated with zero weights. 
    The second half of the statement is equivalent to that the convex envelope of $-f(\cdot)$ over $\conv(\calX)$ is lower semi-continuous. Consider the extended real-valued extension $h:\bbR^n\rightarrow\bbR\cup\{+\infty\}$ of $-f(\cdot)$ to $\bbR^n$, for which $h(\vecx)=-f(\vecx)$ if $\vecx\in\calX$ and $h(\vecx)=+\infty$ otherwise. Note that $h(\cdot)$ is a proper, coercive (by compactness of $\calX$), and lower semi-continuous function, with an epigraph identical to that of $-f(\cdot)$. By \cite{rockafellarVariationalAnalysis1998}[Corollary 3.47], the convex envelope of $h(\cdot)$ is proper, coercive, and lower semi-continuous, and so is the convex envelope of $-f(\cdot)$ over $\conv(\calX)$, which completes the proof.\qed
\end{proof}

\begin{proof}[of Proposition \ref{prop:relax_gap}]
    Concavity and upper semi-continuity of $\sqrt{\hat{q}(\cdot)}$ follow from concavity and monotonicity of $\sqrt{\cdot}$, and concavity and upper semi-continuity of $\hat{q}(\cdot)$. Then $\sqrt{\hat{q}(\cdot)}$ is a concave function majorizing $f(\cdot)=\sqrt{q(\cdot)}$ over $\calX$, and thus majorizing the convave envelope $\hat{f}(\cdot)$ of $f(\cdot)$ over $\conv(\calX)$, which proves the left-hand side inequality of \eqref{ineq:bound}. 
    
    It then remains to prove the right-hand side inequality of \eqref{ineq:bound}.
    Given $\vecx\in\conv(\calX)$, by definition of the concave envelope, there exist $\vecx^1,\ldots,\vecx^K\in\calX$ and $\veclambda^q,\veclambda^f\in\Delta^{K-1}$ such that
    \begin{align}
        \sqrt{\hat{q}(\vecx)}-\hat{f}(\vecx)=\sqrt{\sum_{k=1}^K\lambda^q_k f(\vecx^k)^2}-\sum_{k=1}^K\lambda^f_k f(\vecx^k)\nonumber\\=\sqrt{\sum_{k=1}^K\lambda^q_k f(\vecx^k)^2}-\max_{\lambda\in\Delta^{K-1}}\sum_{k=1}^K\lambda_k f(\vecx^k)\nonumber\\
        \leq \max_{\lambda\in\Delta^{K-1}}\left(\sqrt{\sum_{k=1}^K\lambda_k f(\vecx^k)^2}-\sum_{k=1}^K\lambda_k f(\vecx^k)\right).\label{ineq:opt_bound}
    \end{align}
    Given $\lambda\in\Delta^{K-1}$, let $X$ denote a random variable which takes value $f(\vecx^k)$ with probability $\lambda_k$ for $k=1,\ldots,K$. Let $\mu$ denote its mean and $\sigma^2$ denote its variance, i.e.,
    \begin{align}\label{eq:define mu sigma}
        \mu=\sum_{k=1}^K\lambda_k f(\vecx^k),~\sigma^2={\sum_{k=1}^K\lambda_k f(\vecx^k)^2}-\mu^2.
    \end{align}
    Note that $X$ is supported on $[L,U]$. Then by the Bhatia-Davis inequality \cite{bhatia2000better}, we have\begin{align}\label{ineq:Bhatia-Davis}
        \sigma^2\leq(U-\mu)(\mu-L).
    \end{align}
    Also note that $\mu\in[L,U]$. Combining \eqref{ineq:opt_bound}-\eqref{ineq:Bhatia-Davis}, we have
    \begin{align*}
        \sqrt{\hat{q}(\vecx)}-\hat{f}(\vecx)\leq\max_{\mu\in{[L,U]}}\sqrt{(U-\mu)(\mu-L)+\mu^2}-\mu=\frac{(U-L)^2}{4(L+U)},
    \end{align*}
    where the last equality is obtained by maximizing the middle one-dimensional concave function of $\mu$ over $[L,U]$.\qed
\end{proof}

\section{SOC Reformulation of \eqref{ineq:quad}}\label{apdx:reform_qc}
Since $Q_{yy}$ is positive definite, there exists an invertable matrix $B\in \bbR^{m\times m}$ such that $Q_{yy} = B\trans B$.
Then the inequality can be reformulated as 
\begin{equation*}
    \norm{A\vecx + B\vecy + \vecd}_2 \leq f(x),
\end{equation*}
where
\begin{align*}
    &A = (B\trans)^{-1} Q_{xy}\trans, \quad \vecd = \frac{1}{2}(B\trans)^{-1} \veca_y, \\
    &f(\vecx) = \sqrt{g(\vecx) + \vecx\trans(A\trans A - Q_{xx})\vecx + (2A\trans \vecd -\veca_x)\trans \vecx + \norm{\vecd}_2^2}.
\end{align*}

\section{Instance Generation}\label{apdx:instance_generation}
\subsection{Regular Knapsack}
Following \cite{martelloDynamicProgrammingStrong1999}, we first generate four types of deterministic knapsack instances, which are characterized by different levels of correlation between weights and profits of $N$ items. Within each type, the weight and the profit of each item is independently generated as follows: 
\begin{enumerate}
    \item [1.] Uncorrelated: The weight $w_j$ and profit $p_j$ are integers independently and uniformly distributed in $\{1,2,\ldots,10000\}$;
    \item [2.] Weakly correlated: The weight $w_j$ is an integer uniformly distributed in $\{1,2,\ldots,10000\}$, and the profit is defined as $p_j = \max\{w_j + \kappa_j, 1\}$, where $\kappa_j$ is an integer uniformly distributed in $\{-1000,-999,\ldots,1000\}$;
    \item [3.] Strongly correlated: The weight $w_j$ is an integer uniformly distributed in $\{1,2,\ldots,10000\}$, and $p_j = w_j + 1000$;
    \item [4.] Inverse strongly correlated: The profit $p_j$ is an integer uniformly distributed in $\{1,2,\ldots,10000\}$, and $w_j = p_j + 1000$.
\end{enumerate}
For each instance type and each $N\in \{50, 100, 200\}$, $5$ distinct instances are generated. Specifically, the capacity for the $i$-th instance (where $i = 1,...,5$) is set to $c = \sum_{j=1}^N w_j\cdot i/6$. The weights and profits for the five instances are generated independently, meaning they are not necessarily identical. Following the data generation, all weights and knapsack capacities are divided by $1000$, to enhance the numerical stability.

To adapt these deterministic instances for DRCC-MB-KP, the following modifications are implemented. The last $N/2$ items are designated as continuous items, and their profits are divided by $5$. The mean weight for each item is equal to its deterministic weight, i.e., $\vecmu = \vecw$. The covariance matrix $\Sigma$ is constructed as $\Sigma = U\trans \textrm{diag}(\vecw)\textrm{diag}(\vecw) U/4$, where $U$ is an orthogonal matrix obtained from the decomposition of a randomly generated positive definite matrix. Finally, the knapsack capacity is multiplied by $3/2$, to accommodate a reasonable number of items, especially the discrete items, in the distributionally robust case.

\subsection{Multi-Dimensional Knapsack}
Seven multi-dimensional knapsack instances are adopted from \cite{petersenComputationalExperienceVariants1967} with the number of items $n\in \{6, 10, 15, 20, 28, 39, 50\}$. The number of knapsacks $|\calJ|$ is $10$ for the first five instances and $5$ for the last two. To each of these base instances, $m = \lceil n/2\rceil$ continuous items are introduced. The unit weight for each continuous item is an integer sampled uniformly from $[w_{\min}/2, w_{\max}]$, where $w_{\min}$ (resp. $w_{\max}$) is the minimum (resp. maximum) weight of the discrete items in the base instance. The unit profit of each continuous item is an integer uniformly sampled from $[1, p_{\max}/10]$, where $p_{\max}$ is the maximum profit of discrete items in the base instance. To enhance numerical stability, the weights of all items (both discrete and continuous) are divided by $100$. The mean weight in the stochastic model equals the scaled weight. The covariance matrix is constructed in the same manner as in the regular knapsack case. Finally, the capacities of all knapsacks are divided by $10$. We use $10$ instead of $100$ as the scaling factor to maintain sufficiently large capacities to accommodate a reasonable number of items, as the presence of multiple DRCCs can be quite restrictive.

\section{Experiment Settings}\label{apdx:exp_settings}
All experiments were conducted on a Linux platform (Ubuntu 22.04.5 LTS) equipped with two Intel Xeon Platinum 8575C processors. Gurobi 12.0.1 was used to solve all optimization problems, with a time limit of 3600 seconds and up to 4 threads utilized. We set $\alpha = 0.5\%$ in all the experiments. Trivial initial solutions (setting $\vecx$ and $\vecy$ to $0$) were used when solving {\SOC}, as Gurobi sometimes fails to find feasible solutions for this model.

\end{document}